\documentclass{article}

\usepackage {amsmath}
\usepackage {hyperref}
\usepackage {amsfonts}
\usepackage {amsthm}
\usepackage[flushmargin]{footmisc}
\usepackage[parfill]{parskip}
\usepackage {amssymb}
\usepackage {fullpage}
\usepackage {mathrsfs}
\usepackage {mathtools}
\usepackage {amscd}
\usepackage{bbm}
\usepackage[all,cmtip]{xy}
\DeclareMathAlphabet{\mathpzc}{OT1}{pzc}{m}{it}

\title{On twisted Verlinde formulae for modular categories}
\author{Tanmay Deshpande}
\date{}

\newtheorem {thm} {Theorem} [section]
\newtheorem {prop} [thm] {Proposition}
\newtheorem {conj} [thm] {Conjecture}
\newtheorem {lem} [thm] {Lemma}
\newtheorem {cor} [thm] {Corollary}

\theoremstyle{definition}
\newtheorem {defn} [thm] {Definition}

\newtheorem {prob} [thm]  {Problem}

\newtheorem {rk} [thm]  {Remark}
\newtheorem {ex} [thm] {Example}

\newenvironment{thmbis}[1]
  {\renewcommand{\thethm}{\ref*{#1}$'$}%
   \addtocounter{thm}{-1}%
   \begin{thm}}
  {\end{thm}}

\newcommand{\beq}{\begin{equation}}
\newcommand{\eeq}{\end{equation}}
\newcommand{\bthm}{\begin {thm}}
\newcommand{\ethm}{\end {thm}}
\newcommand{\bprop}{\begin {prop}}
\newcommand{\eprop}{\end {prop}}
\newcommand{\bprob}{\begin {prob}}
\newcommand{\eprob}{\end {prob}}
\newcommand{\bcor}{\begin {cor}}
\newcommand{\ecor}{\end {cor}}
\newcommand{\blem}{\begin{lem}}
\newcommand{\elem}{\end{lem}}
\newcommand{\bdefn}{\begin{defn}}
\newcommand{\edefn}{\end{defn}}
\newcommand{\brk}{\begin{rk}}
\newcommand{\erk}{\end{rk}}

\renewcommand{\subset}{\subseteq}
\newcommand{\xto}{\xrightarrow}

\newcommand{\KabM}{K_{\Qab}(\M)}
\newcommand{\KabC}{K_{\Qab}(\C)}
\newcommand{\KabD}{K_{\Qab}(\D)}

\newcommand{\Kab}{K_{\Qab}}

\renewcommand {\bar} {\overline}
\newcommand{\bpf}{\begin{proof}}
\newcommand{\epf}{\end{proof}}
\newcommand{\bex}{\begin{ex}}
\newcommand{\eex}{\end{ex}}
\newcommand{\rar}[1]{\stackrel{#1}{\longrightarrow}}
\newcommand{\f}{\mathbb}
\newcommand{\fZ}{\mathbb{Z}}
\newcommand{\ZN}{\fZ/N\fZ}

\newcommand{\<}{\langle}
\renewcommand{\>}{\rangle}

\newcommand{\h}{\operatorname}

\newcommand{\un} {\mathbbm{1}}

\renewcommand{\t}{\widetilde}
\renewcommand{\phi} {\varphi}

\newcommand{\M} {\mathscr{M}}

\newcommand{\D} {\mathscr{D}}

\newcommand{\C} {\mathscr{C}}

\renewcommand{\Vec}{\h{Vec}}

\newcommand{\Qab} {\mathbb{Q}^{\h{ab}}}

\newcommand{\A} {\mathcal{A}}

\newcommand{\g} {{\gamma}}

\newcommand{\id}{\operatorname{id}}

\newcommand{\bit}{\begin{itemize}}
\newcommand{\eit}{\end{itemize}}

\newcommand{\Irrep}{\h{Irrep}}

\renewcommand{\O}{\mathcal{O}}

\newcommand{\tr}{\h{tr}}
\newcommand{\Fun}{\h{Fun}}
\newcommand{\Ch}{\h{Ch}}

\newcommand{\bconj}{\begin{conj}}
\newcommand{\econj}{\end{conj}}
\begin{document}
\maketitle
\begin{abstract}
In this note, we describe two analogues of the Verlinde formula for modular categories in a twisted setting. The classical Verlinde formula for a modular category $\C$ describes the fusion coefficients of $\C$ in terms of the corresponding S-matrix $S(\C)$. Now let us suppose that we also have an invertible $\C$-module category $\M$ equipped with a $\C$-module trace. This gives rise to a modular autoequivalence $F:\C\rar{\cong}\C$. In this setting, we can define a crossed S-matrix $S(\C,\M)$. As our first twisted analogue of the Verlinde formula, we will describe the fusion coefficients for $\M$ as a $\C$-module category in terms of the S-matrix $S(\C)$ and the crossed S-matrix $S(\C,\M)$. In this twisted setting, we can also define a twisted fusion $\Qab$-algebra $\Kab(\C,F)$. As another analogue of the Verlinde formula, we describe the fusion coefficients of the twisted fusion algebra in terms of the crossed S-matrix $S(\C,\M)$. 
\end{abstract}

\section{Introduction}\label{s:i}
In this note we will describe two twisted analogues of the Verlinde formula for modular fusion categories, namely Theorems \ref{t:main1} and \ref{t:main2} below. We begin by describing the twisted setting in which we will work throughout this note.

\subsection{Notation and conventions}\label{s:nc}
Let $k$ be an algebraically closed field of characteristic zero. Throughout this note, $\C$ will denote a non-degenerate $k$-linear braided fusion category and $\M$ will denote a $\C$-module category which is invertible when considered as a $\C$-$\C$-bimodule category. We refer to \cite{ENO1}, \cite{ENO2} for the theory of fusion categories and module categories over them. By \cite{ENO2}, the invertible $\C$-module category $\M$ gives rise to an associated braided monoidal autoequivalence $F:\C\rar{\cong}\C$ such that we have functorial crossed braiding isomorphisms
\beq
\beta_{C,M}:C\otimes M \rar{\cong} M\otimes C,
\eeq
\beq
\beta_{M,C}:M\otimes C \rar{\cong} F(C)\otimes M
\eeq
for $C\in \C$, $M\in \M$ satisfying the suitable hexagonal identities. Moreover, as in \cite[\S2.3]{De1}, we can construct a braided $\fZ/N\fZ$-crossed category $\D=\bigoplus\limits_{a\in \fZ/N\fZ}\D_a$ for some positive integer $N$ with trivial component $\D_0\cong \C$ as a braided fusion category and $\D_1\cong \M$ as a module category. We will assume the choice of such a $\D$ throughout the article.

For an abelian category $\A$, let $K(\A)$ denote its Grothendieck group and let $\O_\A$ denote the set of its simple objects up to isomorphism. For an object $A\in \A$, we let $[A]$ denote its image in the Grothendieck group. In our setting, $K(\C)$ is a commutative based ring (see \cite{L}) with basis $\O_\C$ and $K(\M)$ is a $K(\C)$-module with $\fZ$-basis $\O_\M$. For a commutative ring $R$, we let $K_R(\A)=K(\A)\otimes R$.

For most part of this note, we will assume that $\C$ is also equipped with a spherical structure making it a modular fusion category and that $\M$ is equipped with a compatible $\C$-module trace (see \cite{S,De1} for details). In this case, it follows that the braided autoequivalence $F$ is in fact a modular autoequivalence. This extra structure allows us to define categorical dimensions of objects of $\C$ and $\M$. We assume that the $\C$-module trace on $\M$ is normalized according to the conventions of \cite{De1}, so that $\dim \C=\sum\limits_{M\in\O_\M}\dim^2_\M(M)$. With this convention it follows that the categorical dimensions of all objects of $\C$ and $\M$ are totally real cyclotomic integers in $k$. Also, the categorical dimension $\dim \C\in k$ is a totally positive cyclotomic integer. Note that in this setting it is also possible to construct our braided $\fZ/N\fZ$-crossed $\D$ to have a compatible spherical structure (cf. \cite[\S2.4]{De1}).

\subsection{Twisted Verlinde formulae}\label{s:tvf}
We now state the two twisted analogues of the Verlinde formula. The first version below, describes the fusion coefficients for the $\C$-module category $\M$. For simple objects $C\in \O_\C$, $M\in\O_\M$, let
\beq
[C]\cdot [M] = \sum\limits_{N\in\O_\M} a_{C,M}^N [N] \mbox{ in $K(\M)$ with } a_{C,M}^N\in \fZ_{\geq 0}.
\eeq
The first twisted Verlinde formula computes the fusion coefficients $a_{C,M}^N$, i.e. the multiplicity of the simple object $N$ in the product $C\otimes M$.  

Let us assume now that $\C$ is equipped with a spherical structure and that $\M$ is equipped with a normalized $\C$-module trace as in \cite{De1}.  Let $S(\C)$ denote the S-matrix of the modular category $\C$. In this setting, we also have the crossed S-matrix $S(\C,\M)$ defined in \cite{De1}. We will recall the definition and some properties of the S-matrix and the crossed S-matrix in \S\ref{s:csm}. Here we just note that $S(\C)$ is a $\O_\C\times \O_\C$-matrix while $S(\C,\M)$ is a\footnote{Note that the autoequivalence $F:\C\to\C$ induces a permutation $F:\O_\C\to\O_\C$. It can be shown (cf. \cite{De1}) that the number of fixed points $|\O_\C^F|$ equals $|\O_\M|$ although in general there is no canonical bijection between the two sets.} $\O_\C^F\times \O_\M$-matrix and that the matrices $\frac{1}{\sqrt{\dim \C}}\cdot S(\C)$ and $\frac{1}{\sqrt{\dim \C}}\cdot S(\C,\M)$  are unitary (see \cite{De1}). In this setting we have:

\bthm\label{t:main1}
Let $C\in \O_\C$, $M,N\in \O_\M$. Then the multiplicity of $N$ in $C\otimes M$ is given by
\beq\label{e:verlinde1}
a_{C,M}^N=\frac{1}{\dim\C}\sum\limits_{D\in \O_\C^F}\frac{S(\C)_{D,C}\cdot S(\C,\M)_{D,M}\cdot \bar{S(\C,\M)_{D,N}}}{\dim_\C D} 
\eeq
$$
\hspace{27pt}=\frac{1}{\dim\C}\sum\limits_{D\in \O_\C^F}\frac{\bar{S(\C)_{D,C}}\cdot \bar{S(\C,\M)_{D,M}}\cdot {S(\C,\M)_{D,N}}}{\dim_\C D}.
$$
In particular, the two expressions on the right hand side are equal to a non-negative integer.
\ethm
This is the main result of this note. It will be proved in \S\ref{s:fcmc} using results from \cite{De1,De2}. We will also prove a non-spherical version (Theorem \ref{t:main1'}) of this result.

\brk\label{r:crossedsmat}
Recall that the crossed S-matrix $S(\C,\M)$ is only well-defined up to rescaling rows by roots of unity (cf.\cite[Rem. 2.2]{De1}). However it is clear that this does not affect the Verlinde formula (\ref{e:verlinde1}) since any such scaling factors cancel out on the right hand side.
\erk

We now describe another twisted analogue of the Verlinde formula that appeared in \cite{De1}. This formula describes the fusion coefficients of a twisted version of the fusion ring. Let us recall the construction of this twisted fusion ring from \cite{De1}. The autoequivalence $F:\C\rar{}\C$ induces a braided monoidal action of $\fZ$ on $\C$. Consider the Grothendieck ring $K(\C^\fZ)$ of the $\fZ$-equivariantization of $\C$. Now consider the trivial modular category $\Vec$ of finite dimensional $k$-vector spaces and equip it with the trivial action of $\fZ$. We have a ring homomorphism $K(\Vec^\fZ)\rar{}K(\C^\fZ)$. We also have a ring homomorphism $K(\Vec^\fZ)\rar{}k$ that takes the class of the object $(V,\psi:V\cong V)\in \Vec^\fZ$ to $\tr(\psi)\in k$. Define the twisted fusion $k$-algebra $K_k(\C,F)$ to be $K(\C^\fZ)\otimes_{K(\Vec^\fZ)}k$. 

Here is another construction of this algebra. We refer to \cite[\S2.3]{De1} for the details. Recall that for a sufficiently large $N$, we have constructed the braided $\fZ/N\fZ$-crossed category $\D$. In particular, we have a braided action of the finite group $\fZ/N\fZ$ on $\C$ and we can form the braided fusion category $\C^{\fZ/N\fZ}$. Let $\omega$ be a primitive $N$-th root of unity and consider the $\fZ[\omega]$-algebra $K(\C^{\fZ/N\fZ})\otimes_\fZ \fZ[\omega]$. Consider the involution of the ring $\fZ[\omega]$ mapping $\omega$ to $\omega^{-1}$. The rigid duality on $\C^{\fZ/N\fZ}$ provides $K(\C^{\fZ/N\fZ})\otimes_\fZ \fZ[\omega]$ with the structure of a Frobenius $\fZ[\omega]$-$\star$-algebra. Define $K_{\fZ[\omega]}(\C,F)$ to be the quotient of $K(\C^{\fZ/N\fZ})\otimes_\fZ \fZ[\omega]$ by the ideal generated by the element $[(\un,\omega)]-\omega\cdot [\un,\id_\un]$. The twisted fusion $k$-algebra $K_k(\C,F)$ is then defined to be the extension by scalars to $k$. In fact, it will often be more helpful for us to consider the twisted fusion $\Qab$-algebra $\Kab(\C,F):=K_{\fZ[\omega]}(\C,F)\otimes_{\fZ[\omega]}\Qab$, where $\Qab$ is the cyclotomic subfield of $k$ obtained by attaching all roots of unity to $\f{Q}$. Note that $\Qab$ has a canonical field involution denoted by $\bar{(\cdot)}$, mapping a root of unity to its inverse. Note that the rigid duality in $\C^{\fZ/N\fZ}$ equips $\Kab(\C,F)$ with the structure of a Frobenius $\Qab$-$\star$-algebra.

We recall that the twisted fusion algebra $\Kab(\C,F)$ has a basis $\{[C,\psi_C]\}_{C\in \O_C^F}$ parametrized by $\O_C^F$, the set of $F$-stable simple objects of $\C$. For each $C\in \O_\C^F$, $\psi_C:F(C)\rar{\cong}C$ above denotes a (choice of a) $\fZ/N\fZ$-equivariance structure on $C$. This is only well-defined up to scaling by an $N$-th root of unity. Hence, this basis is well-defined only up to scaling by $N$-th roots of unity. 
\brk
From now on, for each $C\in \O_\C^F$ we fix $\psi_C:F(C)\rar{\cong}C$ as above, with $\psi_\un=\id_\un$ and we will work with the basis $\{[C,\psi_C]\}_{C\in \O_C^F}$ of $\Kab(\C,F)$. Note that such a choice is also involved in the definition of the crossed S-matrix $S(\C,\M)$, and we always consider the matrix $S(\C,\M)$ defined using this same  choice of basis. By Remark \ref{r:crossedsmat}, this choice was irrelevant in Theorem \ref{t:main1} above, but will be important below.  
\erk
Our second analogue of the Verlinde formula computes the fusion coefficients for the algebra $\Kab(\C,F)$ in this special basis in terms of the crossed S-matrix $S(\C,\M)$. For $C,C',D\in \O_C^F$, let $a_{C,C'}^D\in \Qab$ be such that
\beq
[C,\psi_C]\cdot[C',\psi_{C'}]=\sum\limits_{D\in \O_\C^F} a_{C,C'}^D[D,\psi_D] \mbox{ in $\Kab(\C,F)$.}
\eeq
From the second construction of the algebra $\Kab(\C,F)$ it is clear that these fusion coefficients $a_{C,C'}^D$ in fact lie in $\fZ[\omega]$. Recall that $S(\C,\M)$ is a $\O_\C^F\times \O_\M$-matrix. We now state the second twisted analogue of the Verlinde formula:
\bthm\label{t:main2} (\cite[Thm. 2.12]{De1}.)
The fusion coefficients for the twisted fusion algebra $\Kab(\C,F)$ are given by
\beq
a_{C,C'}^D=\frac{1}{\dim \C}\sum\limits_{M\in \O_\M}\frac{S(\C,\M)_{C,M}\cdot S(\C,\M)_{C',M}\cdot \bar{S(\C,\M)_{D,M}}}{\dim_\M(M)}
\eeq
for all $C,C',D\in \O_C^F$. In particular, the expression of the right hand side above lies in $\fZ[\omega]$.
\ethm
We will prove a non-spherical analogue of this result in \S\ref{s:fctfa}.

\section{Characters of Grothendieck algebras}\label{s:cga}
In this section we gather general facts about characters of the various Grothendieck algebras and their relationship with S-matrices.

\subsection{Crossed S-matrices}\label{s:csm}
Let us begin by recalling the definition of a crossed S-matrix. We refer to \cite{De1} for details. We work in the setting where $\C$ is a modular category and $\M$ has a normalized $\C$-module trace. In particular, this means that for any endomorphism $\g:M\rar{}M$ in $\M$, the trace $\tr_\M(\g)$ is defined. Let $C\in \O_\C^F$ and let $M\in \O_\M$. Recall that for such a $C$, we have fixed a choice of a $\fZ/N\fZ$-equivariance structure $\psi_C:F(C)\rar{\cong} C$. Consider the composition
\beq
\g_{C,\psi_C,M}: C\otimes M\xto{\beta_{C,M}} M\otimes C\xto{\beta_{M,C}} F(C)\otimes M\xto{\psi_C\otimes \id_M}C\otimes M \mbox{ in } \M.
\eeq
The crossed S-matrix $S(\C,\M)$ is the $\O_\C^F\times \O_\M$-matrix whose $(C,M)$-th entry is
\beq
S(\C,\M)_{C,M}:=\tr_\M(\g_{C,\psi_C,M}),
\eeq
where the trace is computed using the trace in the $\C$-module category $\M$. In the special case when $\M=\C$, we recover the S-matrix $S(\C)$ of $\C$.

We now state some properties of the crossed S-matrices. Recall that we have defined the (commutative) Frobenius $\Qab$-$\star$-algebra (see \cite{A,De1} for details, see also \S\ref{s:fctfa} below) $\Kab(\C,F)$ with our choice of basis $\{[C,\psi_C]|C\in \O^F_\C\}$. On the other hand consider the commutative semisimple $\Qab$-algebra $\Fun_{\Qab}(\O_\M)$ of $\Qab$-valued functions on the set $\O_\M$ with pointwise multiplication. We have:
\bthm\label{t:chartable} (cf. \cite[Thm. 2.9 and 2.12]{De1})
(i) For every $C\in\O_\C^F, M\in \O_\M$, $\frac{S(\C,\M)_{C,M}}{\dim_\C C}$ and $\frac{S(\C,\M)_{C,M}}{\dim_\M M}$ are cyclotomic integers in $k$.\\
(ii) The Fourier transform $\Phi:\Kab(\C,F)\rar{}\Fun_{\Qab}(\O_\M)$, defined on our special basis by 
\beq
\Phi([C,\psi_C]):\O_\M\ni M\mapsto\frac{S(\C,\M)_{C,M}}{\dim_\M M}\in \Qab
\eeq
is an isomorphism of $\Qab$-algebras. In other words, for every $M\in \O_\M$, $[C,\psi_C]\mapsto \frac{S(\C,\M)_{C,M}}{\dim_\M M}$ defines a character $\phi_M:\Kab(\C,F)\rar{} \Qab$ and this defines an identification $\O_\M\cong\Irrep(\Kab(\C,F))$.\\
(iii) The matrix $\frac{1}{\sqrt{\dim \C}}\cdot {S(\C,\M)}$ is unitary:
\beq
S(\C,\M)\cdot \bar{S(\C,\M)}^T = \dim \C\cdot I = \bar{S(\C,\M)}^T\cdot {S(\C,\M)}.
\eeq
\ethm
The above result says that the crossed S-matrix is essentially the character table of the twisted fusion algebra $\Kab(\C,F)$ in the basis $\{[C,\psi_C]\}_{C\in \O_\C^F}$. Next, we relate the crossed S-matrix to a certain twisted character table.

\subsection{Twisted characters and crossed S-matrices}\label{s:tccsm}
We have constructed an auxiliary braided $\fZ/N\fZ$-crossed category $\D$. Consider the $\fZ/N\fZ$-graded based $\Qab$-algebra $\Kab(\D)=\bigoplus\limits_{a\in \ZN}\Kab(\D_a)$ with basis $\{[D]\}_{D\in\O_\D}$. It is known (see \cite[Cor. 8.53]{ENO1}) that all the irreducible representations of $\KabD$ and $\Kab(\C,F)$ are in fact defined over $\Qab$ (see \cite{De1}). 

Now $\Kab(\D)$ is a $\fZ/N\fZ$-graded Frobenius $\Qab$-$\star$-algebra and hence we have the notion of twisted characters and their orthogonality relations (cf. \cite{A}) which we recall now. By \cite{Da,De2}, the $\fZ/N\fZ$-graded algebra $\Kab(\D)$ induces a partial $\ZN$-action on the set of irreducible representations of the identity component $\Kab(\C)$. On the other hand, the braided monoidal $\ZN$ action on $\C$ induces an action of $\ZN$ on the algebra $\Kab(\C)$ and hence also on $\Irrep(\Kab(\C))$. However, this action does not agree with the partial action obtained above. Using the argument from \cite[Lem. 3.5]{De2} we can describe the partial action:
\blem\label{l:partialaction}
For $a\in \ZN, \rho \in \Irrep(\Kab(\C))$, let ${ }^{(a)}\rho\in \Irrep(\KabC)\cup \{0\}$ denote the partial action. For $a\in\ZN$, we have the algebra automorphism $F^a:\KabC\rar{}\KabC$. Then we have
\beq
{ }^{(a)}\rho = 
\begin{dcases*}
\rho & if $\rho \circ F^a=\rho$,\\
0 & else.
\end{dcases*}
\eeq
In particular, the fixed points of the partial action and the action match, so the notation $\Irrep(\KabC)^{\ZN}$ is unambiguous.
\elem
By the results of $\cite{Da,De2}$, each $\rho\in \Irrep(\KabC)^{\ZN}$ can be extended to a 1-dimensional character $\t\rho:\KabD\rar{}\Qab$. There are in fact $N$ such extensions and these extensions differ on $\Kab(\M=\D_1)$ up to scaling by the $N$-th roots of unity. For $\rho\in \Irrep(\KabC)^{\ZN}$, its twisted character is the linear functional $\t\chi_\rho:=\t\rho|_{\KabM}:\KabM\rar{}\Qab$. As we have noted, the twisted character $\t\chi_\rho$ is well-defined up to scaling by $N$-th roots of unity. 

Since $\KabD$ is a Frobenius $\Qab$-$\star$-algebra, we will identify (as in \cite{De2}) the linear dual $\KabM^*$ with $\Kab(\M^{-1}=\D_{-1})$. The basis $\{[D]\}_{D\in \O_\D}$ of $\KabD$ is orthonormal with respect to the standard positive definite Hermitian form denoted by $\<\cdot,\cdot\>$ (cf. \cite{De2}). Hence, under the identification $\KabM^*=\Kab(\M^{-1})$, a functional $\phi\in \KabM^*$ corresponds to the element $\sum\limits_{M\in\O_\M}\phi([M])[M^*]\in \Kab(\M^{-1})$. 

The twisted characters $\{\t\chi_\rho|\rho\in\Irrep(\KabC)^{\ZN}\}$ form an orthogonal basis of $\KabM^*=\Kab(\M^{-1})$ with respect to the standard positive definite Hermitian form (see \cite{A,De2} for details). For each $\rho\in\Irrep(\KabC)^{\ZN}$, let $\t\alpha_\rho\in \Kab(\M^{-1})$ be the element corresponding to the twisted character $\t\chi_\rho\in \KabM^*$, namely
\beq\label{e:talpha}
\t\alpha_\rho=\sum\limits_{M\in \O_\M} \t\chi_\rho([M])[M^*]\in \Kab(\M^{-1}).
\eeq 
Similarly, for each $\rho\in \Irrep(\KabC)$, we have the element $\alpha_\rho=\sum\limits_{C\in \O_\C} \rho([C])[C^*]\in \KabC$ corresponding to $\rho\in \KabC^*$. Note that $\Kab(\M^{-1})$ is a $\KabC$-module. In this notation, we have the following:
\bthm\label{t:idempotents} (cf. \cite{L,O},\cite[\S2.7]{De2}.)
(i) Let $\rho,\rho'\in\Irrep(\KabC)$, then $\rho'(\alpha_\rho)=0$ if $\rho'\neq\rho$ and $f_\rho:=\rho(\alpha_\rho)\in\Qab$ is a totally positive cyclotomic integer known as the formal codegree of $\rho$. In other words, $\frac{\alpha_\rho}{f_\rho}\in \KabC$ is the minimal idempotent corresponding to $\rho\in\Irrep(\KabC)$.\\
(ii) For $\rho,\rho'\in\Irrep(\KabC)$, we have orthogonality of characters, $\<\alpha_\rho,\alpha_{\rho'}\>=\delta_{\rho\rho'}\cdot f_\rho$.\\
(iii) Let $\rho\in \Irrep(\KabC)^{\ZN}$ and $\t\alpha_{\rho}\in \Kab(\M^{-1})$ corresponding to the twisted character $\t\chi_\rho\in\KabM^*$. Then for any $\rho'\in\Irrep(\KabC)$, we have 
\beq
\alpha_{\rho'}\cdot \t\alpha_\rho =
\begin{dcases*}
0 & if $\rho'\neq \rho$,\\
f_\rho\cdot \t\alpha_\rho & if $\rho'=\rho$.
\end{dcases*}
\eeq
We have a direct sum decomposition
\beq
\Kab(\M^{-1})=\bigoplus\limits_{\rho\in\Irrep(\KabC)^{\ZN}} \Qab\cdot \t\alpha_\rho \mbox{ as a $\KabC$-module.}
\eeq
(iv) For $\rho,\rho'\in\Irrep(\KabC)^{\ZN}$, we have orthogonality of twisted characters, $\<\t\alpha_\rho,\t\alpha_{\rho'}\>=\delta_{\rho\rho'}\cdot f_\rho$.
\ethm

Now assume that $\C$ is modular, $\M$ has a $\C$-module trace and that $\D$ is chosen to have a compatible spherical structure. In this case, we can identify $\O_\C\cong \Irrep(\KabC)$ (cf. \cite{ENO1}, also Thm. \ref{t:chartable} above) and the action of $\ZN$ on both sides matches. For $C\in\O_\C$, the map $\phi_C:\KabC\rar{}\Qab, [D]\mapsto \frac{S(\C)_{C,D}}{\dim_\C C}$ defines the corresponding irreducible representation. For each $C\in \O_\C^F\cong\Irrep(\KabC)^{\ZN}$, we have chosen an equivariantization $(C,\psi_C)\in \C^{\ZN}\subset \D^{\ZN}$. As observed in \cite[\S4.4]{De2}, this (choice of equivariantization) determines an extension of $\phi_C$ to a character $\t\phi_C=\phi_{C,\psi_C}:\KabD\rar{}\Qab.$ As in loc. cit., we obtain:
\bthm\label{t:twistedchartable}
(i) For each $C\in \O_\C^F\cong \Irrep(\KabC)^{\ZN}$, the twisted character $\t\chi_{C}=\t\chi_{\phi_C}:\KabM\to\Qab$ (determined by our fixed choice of $\psi_C:F(C)\xto{\cong} C$) is given by
\beq
\t\chi_C([M])=\frac{S(\C,\M)_{C,M}}{\dim_\C C} \mbox{ for each } M\in \O_\M.
\eeq
Hence in the spherical setting, we have $\t\alpha_C=\sum\limits_{M\in \O_\M}\frac{S(\C,\M)_{C,M}}{\dim_\C C}\cdot [M^*]\in \Kab(M^{-1})$.\\
(ii) For any $C\in\O_\C\cong \Irrep(\KabC)$, the formal codegree $f_C=f_{\phi_C}$ of the representation $\phi_C$ equals the totally positive cyclotomic integer $\frac{\dim \C}{\dim^2_\C C}$.
\ethm
\bpf
The proof of (i) follows from the argument in \cite[\S4.3]{De2}. To prove (ii), take $\M=\C$. Then by Theorems \ref{t:chartable}(iii), \ref{t:idempotents}(ii) and statement (i) above, we obtain $f_C=\<\alpha_C,\alpha_C\>=\sum\limits_{C'\in\O_\C}\frac{S(\C)_{C,C'}\cdot \bar{S(\C)_{C,C'}}}{\dim_\C^2 C}=\frac{\dim \C}{\dim_\C^2 C}$.
\epf
The above result says that in the spherical case, the crossed S-matrix is essentialy the twisted character table of $\KabC$. 

\section{Proof of the main results}
We can now complete the proof of the main Theorem \ref{t:main1}. Theorem \ref{t:main2} is already proved in \cite{De1}. Here, we will state and prove  versions of both the twisted Verlinde formulae without assuming the existence of spherical structures.

\subsection{Fusion coefficients for the module category}\label{s:fcmc}
Let us first prove a version of Theorem \ref{t:main1} without spherical structures. In this case the S-matrix and crossed S-matrix are not available. We will instead work with the character table and the twisted character table of $\KabC$. 

Let us now consider the two character table matrices. Let $\Ch$ be the $\Irrep(\KabC)\times \O_\C$ matrix with $(\rho,C)$-th entry defined as $\Ch_{\rho,C}=\rho([C])$. Similarly, let $\t\Ch$ be the $\Irrep(\KabC)^{\ZN}\times \O_\M$ matrix with $(\rho,M)$-th entry $\t\Ch_{\rho,M}:=\t\chi_{\rho}([M])$. We will now describe the product $[C^*]\cdot[M^*]\in \Kab(\M^{-1})$ for $C\in \O_\C, M\in \O_\M$. 

By Theorem \ref{t:idempotents}(ii) and (iv), we have $\Ch\cdot \bar{\Ch}^T=\h{Codeg}$ and $\t\Ch\cdot \bar{\t\Ch}^T=\t{\h{Codeg}}$, where $\h{Codeg}$ (resp. $\t{\h{Codeg}}$) is the $\Irrep(\KabC)\times \Irrep(\KabC)$ (resp. $\Irrep(\KabC)^{\ZN}\times \Irrep(\KabC)^{\ZN}$) diagonal matrix with $(\rho,\rho)$-th entry $f_\rho$, the formal codegree of $\rho$.

By our definitions above (see also \ref{e:talpha}) we have $\left(\begin{array}{c}\vdots\\\ \alpha_\rho\\\vdots\end{array}\right)=\Ch\cdot\left(\begin{array}{c}\vdots\\\ [C^*]\\\vdots\end{array}\right)$ and $\left(\begin{array}{c}\vdots\\\ \t\alpha_\rho\\\vdots\end{array}\right)=\t\Ch\cdot\left(\begin{array}{c}\vdots\\\ [M^*]\\\vdots\end{array}\right)$. Hence $\bar{\Ch}^T\cdot\h{Codeg}^{-1}\cdot\left(\begin{array}{c}\vdots\\\ \alpha_\rho\\\vdots\end{array}\right)=\bar{\Ch}^T\cdot\left(\begin{array}{c}\vdots\\\ \frac{\alpha_\rho}{f_\rho}\\\vdots\end{array}\right)=\left(\begin{array}{c}\vdots\\\ [C^*]\\\vdots\end{array}\right)$ and $\bar{\t\Ch}^T\cdot\t{\h{Codeg}}^{-1}\cdot\left(\begin{array}{c}\vdots\\\ \t\alpha_\rho\\\vdots\end{array}\right)=\left(\begin{array}{c}\vdots\\\ [M^*]\\\vdots\end{array}\right)$.
We can now prove the following non-spherical analogue of Theorem \ref{t:main1}
\begin{thmbis}{t:main1}\label{t:main1'}
For $C\in \O_\C$ and $M,N\in \O_\M$, we have
\beq
a_{C,M}^N=a^{N^*}_{C^*,M^*}=\sum\limits_{\substack{\rho\in \\\Irrep(\KabC)^{\ZN}}}\frac{\bar{\rho([C])}\cdot\bar{\t\chi_\rho([M])}\cdot\t\chi_\rho([N])}{f_\rho}=\sum\limits_{\substack{\rho\in \\\Irrep(\KabC)^{\ZN}}}\frac{{\rho([C])}\cdot{\t\chi_\rho([M])}\cdot\bar{\t\chi_\rho([N])}}{f_\rho}.
\eeq
\end{thmbis}
\bpf
For $C\in\O_\C,M\in \O_\M$, let us first compute $[C^*]\cdot[M^*] \in \Kab(\M^{-1})$. By the previous equations
\beq
[C^*]\cdot[M^*]=\left(\sum\limits_{\rho\in \Irrep(\KabC)}\bar{\rho([C])}\cdot\frac{\alpha_\rho}{f_\rho}\right)\cdot\left(\sum\limits_{\rho\in \Irrep(\KabC)^{\ZN}}\bar{\t\chi_\rho([M])}\cdot\frac{\t\alpha_\rho}{f_\rho}\right)
\eeq
\beq
=\sum\limits_{\rho\in \Irrep(\KabC)^{\ZN}}\frac{\bar{\rho([C])}\cdot\bar{\t\chi_\rho([M])}}{f_\rho}\cdot{\t\alpha_\rho} \mbox{\hspace{20pt}}\mbox{ $\cdots$ by Thm. \ref{t:idempotents}(iii)}
\eeq
\beq
=\sum\limits_{\substack{{\rho\in \Irrep(\KabC)^{\ZN}}\\{N\in \O_\M}}}\frac{\bar{\rho([C])}\cdot\bar{\t\chi_\rho([M])}}{f_\rho}\cdot\t\chi_\rho([N])\cdot [N^*]
\eeq
Hence for each $N\in \O_\M$, the multiplicity of $[N^*]$ in $[C^*]\cdot[M^*]$ is given by 
\beq
a^{N^*}_{C^*,M^*}=\sum\limits_{{\rho\in \Irrep(\KabC)^{\ZN}}}\frac{\bar{\rho([C])}\cdot\bar{\t\chi_\rho([M])}\cdot\t\chi_\rho([N])}{f_\rho}.
\eeq
Now it is clear that we have $(C\otimes M)^*\cong C^*\otimes M^*$ as objects of $\M^{-1}$. Hence we have $a_{C,M}^N=a^{N^*}_{C^*,M^*}$. Moreover, since the multiplicities are integers, we can take the `complex conjugation' in the summation and (noting that $f_\rho$ is a totally positive cyclotomic integer) the theorem follows.
\epf

\begin{proof}[Proof of Theorem \ref{t:main1}]
Now we assume that we are in the spherical setting. In this setting we have identifications $\O_\C\cong \Irrep(\KabC)$ and $\O_\C^F\cong \Irrep(\KabC)^{\ZN}$. Now suppose that $D\in \O_\C^F$ corresponds to $\rho\in \Irrep(\KabC)^{\ZN}$, then we have 
\beq
\rho([C])=\frac{S(\C)_{D,C}}{\dim_\C D}, \mbox{ }\t\chi_\rho([M])=\frac{S(\C,\M)_{D,M}}{\dim_\C D}, \mbox{ } \t\chi_\rho([N])=\frac{S(\C,\M)_{D,N}}{\dim_\C D} \mbox{ and } f_\rho = \frac{\dim \C}{\dim^2_\C D}.
\eeq
The first equality above follows from either Theorem \ref{t:chartable}(ii) or \ref{t:twistedchartable}(i) applied to the case $\M=\C$. The next three equalities follow from Theorem \ref{t:twistedchartable}(i) and (ii). This combined with Theorem \ref{t:main1'} completes the proof.
\epf

\subsection{Fusion coefficients for the twisted fusion algebra}\label{s:fctfa}
Theorem \ref{t:main2} is proved in \cite{De1}. Here we will only derive a non-spherical version of this Verlinde formula. The role of the S-matrix will be played by the character table of the twisted fusion algebra $\Kab(\C,F)$ in the basis $\{[C,\psi_C]\}_{C\in\O_\C^F}$. 

We still have that $\Kab(\C,F)$ is a commutative Frobenius $\Qab$-$\star$-algebra. Let us recall what this means. We refer to \cite{De1,A} for details. Firstly, we have a non-degenerate linear functional $\lambda:\Kab(\C,F)\rar{}\Qab$ defined by $\lambda([\un,\id_\un])=1$ and $\lambda([C,\psi_C])=0$ for $\un\ncong C\in \O_\C^F$ providing $\Kab(\C,F)$ with the structure of a Frobenius algebra. In particular we have an identification $\Kab(\C,F)\cong \Kab(\C,F)^*$ as $\Kab(\C,F)$-modules. The two bases $\{[C,\psi_C]\}_{C\in \O_\C^F}$ and $\{[C^*,{\psi^*_C}^{-1}]\}_{C\in \O_\C^F}$ of the Frobenius algebra $\Kab(\C,F)$ are dual in the usual sense. Moreover, there is a $\Qab$-semilinear anti-involution $(\cdot)^*:\Kab(\C,F)\rar{}\Kab(\C,F)$ mapping a basis element $[C,\psi_C]\mapsto [C^*, {\psi_C^*}^{-1}]$ and extended semilinearly. Then the Hermitian form $\<\cdot,\cdot\>$ on $\Kab(\C,F)$ defined by $\<a,b\>=\lambda(ab^*)$ is (totally) positive definite and our special basis $\{[C,\psi_C]\}_{C\in \O_\C^F}$ is orthonormal.

We reiterate that by \cite{ENO1}, all irreducible characters of $\Kab(\C,F)$ are defined over $\Qab$. For a character $\phi:\Kab(\C,F)\rar{}\Qab$, we have the corresponding element $\alpha_\phi=\sum\limits_{C\in \O_\C^F}{\phi([C^*,{\psi_C^*}^{-1}])}[C,\psi_C]=\sum\limits_{C\in \O_\C^F}{\phi([C,\psi_C])}[C^*,{\psi_C^*}^{-1}]$ in $\Kab(\C,F)$. From generalities about Frobenius $\star$-algebras and the definition of $\Kab(\C,F)$, it follows that for $\phi' \in \Irrep(\Kab(\C,F))$, $\phi'(\alpha_\phi)=0$ if $\phi'\neq \phi$ and $f_\phi:=\phi(\alpha_\phi)$ is a totally positive cyclotomic integer known as the formal codegree of $\phi$. We have the orthogonality relations $\<\alpha_\phi,\alpha_{\phi'}\>=\delta_{\phi\phi'}\cdot f_\phi$. The element $\frac{\alpha_\phi}{f_\phi}\in \Kab(\C,F)$ is the minimal idempotent corresponding to the irreducible representation $\phi$. It follows that (see \cite[Lem. 2.5]{A}) $\alpha_\phi^*=\alpha_\phi$ and hence $\phi([C^*,{\psi_C^*}^{-1}])=\bar{\phi([C,\psi_C])}$.

Let $\Ch_F$ be the $\Irrep(\Kab(\C,F))\times \O_\C^F$ matrix with $(\phi,C)$-th entry $\phi([C,\psi_C])$. We have $\Ch_F\cdot\bar{\Ch_F}^T=\h{Codeg}_F$, where $\h{Codeg}_F$ is the $\Irrep(\Kab(\C,F))\times \Irrep(\Kab(\C,F))$ matrix whose $(\phi,\phi)$-th entry is the formal codegree $f_\phi$. By our definitions $\left(\begin{array}{c}\vdots\\\ \alpha_\phi\\\vdots\end{array}\right)=\Ch_F\cdot\left(\begin{array}{c}\vdots\\\ [C^*,{\psi^*_C}^{-1}]\\\vdots\end{array}\right)=\bar{\Ch_F}\cdot\left(\begin{array}{c}\vdots\\\ [C,{\psi_C}]\\\vdots\end{array}\right)$. Hence  $\Ch_F^T\cdot\left(\begin{array}{c}\vdots\\\ \frac{\alpha_\phi}{f_\phi}\\\vdots\end{array}\right)=\left(\begin{array}{c}\vdots\\\ [C,{\psi_C}]\\\vdots\end{array}\right)$. Hence proceeding as in the proof of Theorem \ref{t:main1'} we obtain (compare with \ref{t:main2}):
\begin{thmbis}{t:main2}
For $C,C',D\in \O_C^F$ we have
\beq
a_{C,C'}^D=\sum\limits_{\phi\in \Irrep(\Kab(\C,F))}\frac{\phi([C,\psi_C])\cdot \phi([C',\psi_{C'}])\cdot\bar{\phi([D,\psi_D])}}{f_\phi}.
\eeq
\end{thmbis}
\brk
Theorem \ref{t:main2} (in the spherical setting) follows from this using Theorem \ref{t:chartable}(ii) and its consequence that for $M\in \O_\M$, the formal codegree of the associated representation $\phi_M$ is given by $f_{\phi_M}=\frac{\dim\C}{\dim_\M^2 M}$.
\erk

\end{document}